   \input amstex
\documentstyle{amsppt}
\magnification=\magstep1
\parindent=1em
\baselineskip=18pt
\CenteredTagsOnSplits
\NoBlackBoxes
\nopagenumbers
\NoRunningHeads
\pageno=1
\footline={\hss\tenrm\folio\hss}

\topmatter
\title {Approximation properties $ \operatorname{\bold{AP_s}}$ and
            ${\ssize\bold p}$-nuclear operators
     \eightpoint   \bf (the case where $ 0<{\ssize\bold s}\le1$) }
\endtitle
\author { O. I. Reinov${{ }^\dag}$}  \endauthor
\address\newline
Oleg I. Reinov \newline
Department of Mathematics\newline
St Petersburg University\newline
St Peterhof, Bibliotech pl 2\newline
198904  St Petersburg, Russia
\endaddress

\email
orein\@orein.usr.pu.ru
\endemail

\thanks
${{ }^\dag}$ This work was done with partial support of the Ministry of the
general and professional education of Russia (Grant 97-0-1.7-36) and
FCP ``Integracija", reg. No. 326.53.
\endthanks

\abstract
Among other things, it is shown that there exist Banach spaces
$Z$ and $W$ such that $Z^{**}$ and $W$ have bases, and for every
$p\in[1,2)$ there is an operator $T:W\to Z$ that is not $p$-nuclear
but $T^{**}$ is $p$-nuclear.
\endabstract
    \endtopmatter

\document
\baselineskip=18pt
\footnote""{${ }^\ddag$
}
\footnote""{${ }$
}

\magnification=\magstep1

 \def\a{\alpha}                         \def\ot{\otimes}
 \def\la{\lambda}                       \def\wh{\widehat}
 \def\ffi{\varphi}                      \def\wt{\widetilde}
 
\def\Q{\quad\blacksquare}
\def\A{\operatorname{AP}}
\def\AP#1{\operatornamewithlimits{AP_#1}}
\def\ap[#1#2]{\operatornamewithlimits{AP}_{#1,#2}}

\def\nr[#1#2]{\operatornamewithlimits{N^{reg}_{#1,#2}\,}}
\def\NR[#1]{\operatornamewithlimits{N^{reg}_{#1}\,}}
\def\n#1{\operatornamewithlimits{N_#1}\,}
\def\nd[#1#2]{\operatornamewithlimits{N^{dual}}_{#1,#2}\,}
\def\ND[#1]{\operatornamewithlimits{N^{dual}}_{#1}\,}
\def\nu[#1]{\operatornamewithlimits{N}^#1\,}

\def\({\left(}\def\){\right)}\def\[{\left[}\def\]{\right]}
\def\tr{\operatorname{trace}\,}

\bigskip

We hold standard notation of the geometrical theory of operators
in Banach spaces. The classical reference book on the theory of operator
ideals is the A. Pietsch monograph [7].
Notations and terminology we use can be found, for instance, in
[7], [8], [10], [11].
For our purposes, it is enough only to recall that if $ X, Y$ is a pair
of Banach spaces and $ p>0,$ then $ N_p(X, Y)$ denotes
the space of all $ p$-nuclear operators from $ X$ to $ Y,$ and
$ X^*\wh\ot_p Y$ denotes the associated with it $ p$-projective tensor
product. And one more important reminder. If $ J$ is an operator ideal,
then $ J^{ \operatorname{ reg}}(X, Y)$ denotes
the space of all operators $ T$ from $ X$ into $ Y,$ for which
$ \pi_Y\,T\in J(X, Y^{**}),$ where $ \pi_Y$ is
the canonical isometric imbedding of the space $ Y$
into its second dual $ Y^{**}.$
\smallpagebreak

In this note we are interesting mainly in the following questions:
\roster
\item "{1)}"
Under which conditions on Banach spaces $\,X,Y\,$
the canonical mapping
$X^*\widehat \otimes _p \,Y\to\operatornamewithlimits{N_p}\left(X,Y \right)$
is one--to--one?
\item "{2)}"
Under which conditions on Banach spaces $\,X,Y\,$
and on positive numbers $\,s,r\,$ \
$\operatornamewithlimits{N_s^{reg}}\left(X,Y \right)\subset
\operatornamewithlimits{N_r}\left(X,Y \right)$?
\endroster
Both questions were investigated rather explicitly, e.g., in papers
[10], [11].
The question 2) was considered also in [2] for $r=s=1$ and in [8]
for $r=s>1$.
For some other values of parameters $p,r,s$ both questions were
analyzed in the paper [3]. So, the well known A. Grothendieck's theorem
on 2/3 asserts that always one has $ \operatorname{ N_{2/3}^{reg}}\subset
\operatorname{ N}_1$ (somewhat below a proof of this assertion can be seen
in the remark 2).

As more examples let us note also the following facts:
\roster
\item "{a)}"
if $X,Y$ are arbitrary Banach spaces then the canonical mapping
$X^*\widehat \otimes _{2/3} \,Y \to
\operatorname{L}\left(X,Y \right)$  is always one--to--one (see [3], [9],
and also the corollary $ 1'$ below);
\item "{b)}"
if every finite--dimensional subspace $E$ of a space $Y\ $ is
 $\ C\left({\operatorname{dim}E}\right)^\alpha $-complemented in $ Y$
(where $0<\alpha\le 1/2$), then for $s,$ $\ 1/s=1+\alpha,$ and
for every Banach space $X$ the canonical mapping
$X^*\widehat \otimes _{s} \,Y \to\operatorname{L}\left(X,Y \right)$
is one--to--one (see [10]);
\item "{c)}"
if $0<s\le1$ and $1/s=1/2+1/r,$  then
$\operatornamewithlimits{N_s^{reg}}\left(X,Y \right)\subset
\operatornamewithlimits{N_r}\left(X,Y \right)$
for any Banach spaces $X,Y$ (see [3]);
\item "{d)}"
if $p\ge 1$ and $X^*\in \operatorname{AP},$ then
$\operatornamewithlimits{N_p^{reg}}\left(X,Y \right)\subset
\operatornamewithlimits{N_p}\left(X,Y \right)$
for each Banach space $Y$ (see [3] concerning the case
$p=1;$ what about the case $p>1$, here the proof is carried out
by the same scheme that was in [3], so we are not going to repeat it).
\endroster
\remark {\bf Remark {\rm 1}}
A similar assertion as d) with assumption that $ Y^{**}\in
\operatorname{ AP}$   (instead of $ X^*$) was formulated, without a proof,
in [3] for
$ p=1$ and in [5], with a proof based on [3],
for $ p>1.$ Below we will see that these last facts
are not valid at least for the cases when $ 1\le p<2$
(however, concerning the case $ p=1$ see the note [6]).
\endremark
In [10], among others, it is shown that the assertions \, (a)--(c) can not
be improved ``in the scale of the spaces $ l_p.$" Below we shall refine
the assertions
(a)--(c) for ``the scale of Lorentz spaces"\footnote{
But not in the very full generality;
below we shall be interesting in addition only in spaces
$ l_{s,\infty}$ for $ s< 1;$ but already in this case
there are non entirely evident generalizations;
see, e.g., theorem 4 below.},
and shall show that the new obtained results are unimprovable
in this cases too.

We shall see also (and this is, perhaps, the main thing in the paper), that
all the ``negative" results of [8], [10], corresponding to the assertions
(a)--(d), will take a place in the stronger versions. The first example of
such a kind of the ``sharpened" counterexamples
may be found in the paper [6],
where it was shown that there exists a Banach space $ Z$ with special
properties: (i) $Z^{**}$ is separable and has a basis
(and $ Z^{***}$ does not possess the approximation property);
(ii)
$\operatornamewithlimits{N_1^{reg}}\left(Z^{**},Z \right)\not\subset
\operatornamewithlimits{N_1}\left(Z^{**},Z \right)$.
This is an example demonstrating that the second part of a Grothendieck's
assertion [3] (ch.1, p. 86, proposition 15) is not valid.
Below we will extend the mentioned counterexample to the case of the
exponents of nuclearity $ p,$ bigger than one and less than two
(but, certainly, differed from two; the corresponding counterexamples
for $ p>2$ are the theme of a forthcoming paper); see the theorem 5
and its corollary 3.

Below we often shall use two following facts
obtained in the paper [4], --- the first one is exactly
the theorem 1 in [4], and the second one is the corollary 1 from there
(the existence of a basis in the space $ Z^{**}$
in the lemma 2 follows from the definition of the space $ Z$
in the proof of that corollary in [4]).

\proclaim {\bf Lemma 1}\it
For each separable Banach space $ X$
there exist a separable Banach space $ E$
and a linear homomorphism $ \psi: E^*\to X$ such that
$ E^*$ has a basis and $ E^{**}=\psi^*(X^*)\oplus \pi_E(E).$
\endproclaim\rm

\proclaim {\bf Lemma 2}\it
For each separable Banach space $ Y$
there exist a separable Banach space $ Z$
and a linear homomorphism  $ \psi: Z^{**}\to Y$ such that
$ Z^{**}$ has a basis
{\rm(}and, consequently, is separable{\rm)},
the kernel of the homomorphism $ \psi$ is equal to $ \pi_Z(Z)$
and $ \psi ^*(Y^*)$ is complemented in $ Z^{***}.$

\endproclaim\rm

\proclaim {\bf Corollary 1}\it
For each separable Banach space $ Y$
there exists a separable Banach space $ Z$
such that $ Z^{**}$ has a basis , $ Z^{**}/\pi_Z(Z)$
is isomorphic to the space $ Y,$\ $ \pi_Z(Z)^{\perp}$
is complemented in $ Z^{***}$
and the following conditions are fulfilled:\newline
\phantom{vvv}{\rm(}a{\rm)}
for every $ p\in (0,1]$\ and every Banach space
$ E$ the tensor product $E^*\wh\ot_p Y$ is isomorphic to the factor space
$E^*\wh\ot_p \(Z^{**}/\pi_Z(Z)\),$ which is, in turn, is isomorphic to
$$ N_p\(E, Z^{**}\)/\[\NR[p]\(E,\pi_Z(Z)\)
\bigcap  N_1(E, \pi_Z(Z)) \];$$
\phantom{vvv}{\rm(}b{\rm)}
for every $ p\in (0,1]$\ and every Banach space
$ E$ the space $ N_p(E,Y)$ is isomorphic to the space
$ N_p\(E,Z^{**}/\pi_Z(Z)\),$ which is, in turn, is isomorphic
to the factor space $ N_p\(E, Z^{**}\)/\NR[p]\(E,\pi_Z(Z)\).$

In addition, all the isomorphisms in the assertions {\rm(}a{\rm)}
and {\rm(}b{\rm)}
are ``canonical", that is, are generated by the given homomorphism
$ Z^{**}\to Z^{**}/\pi_Z(Z)\to Y$
{\rm (details are in the proof)}.
\endproclaim\rm
\demo{\it Proof}
Let $ Z$ be the space from lemma 2. Denote, for the simplicity,
$ Z^{**}/\pi_Z(Z)$ by $ F.$

Quotient map $ \psi: Z^{**}\to Z^{**}/\pi_Z(Z)=F $
induces naturally the map
$ \Psi: N_p(E, Z^{**})\to E^*\wh\ot_p F,\,$ $ \Psi(T):= \psi\circ T.$
Since $ p\le 1,$ it is clear that this map $ \Psi$ is a surjection.
Let us see what is the kernel of this map.

Let $ U\in \operatorname{ Ker}\,\Psi.$ This means that for every
$ A\in \operatorname{ L}(F, E^{**})$ \, $ \tr A\circ\psi\circ U$
is equal to zero; or,
that for each $ B\in \operatorname{ L}(Z^{**}, E^{**})$
such that $ B|_{\pi_Z(Z)}=0,$ \, $ \tr B\circ U$ is zero.
It follows that $ U(E)\subset \pi_Z(Z)$  (in other case
one can find an operator
$ R\in \operatorname{ L}(Z^{**}, E^{**})$ such that $ R|\pi_Z(Z)=0$ and
$ \tr R\circ U\neq 0$).

So, $ U\in \NR[p](E, \pi_Z(Z)).$  We have to refine this inclusion:
let us prove that in fact
$ U\in N_1(E, \pi_Z(Z)).$
Suppose that $ U\notin N_1(E, \pi_Z(Z)).$ Then one can find
an operator $ B\in \operatorname{ L}(Z^{**}, E^{**}),$ such that
$ \tr B\circ U =1,$ but $ \tr B\circ T =0$ for every
$ T\in N_1(E, \pi_Z(Z)).$ Since $ \pi_Z(Z)$ has the Grothendieck
approximation property, it follows now that
$ B|_{\pi_Z(Z)}=0,$ what is contradicted to the choice and properties
of $ U.$
So, $ U\in N_1(E, \pi_Z(Z)).$ On the other hand, evidently, for every
$ V\in N_1(E, \pi_Z(Z))$ and $ B\in \operatorname{ L}(Z^{**}, E^{**}),$
where $ B|_{\pi_Z(Z)}=0,$ we have: $ \tr B\circ V=0.$
Thus $ \operatorname{ Ker}\,\Psi =
\NR[p]\(E,\pi_Z(Z)\) \cap \( N_1(E, \pi_Z(Z))\).$ It follows from here
the assertion (a).

Now consider the quotient map
$$ \Psi_0:\qquad N_p(E, Z^{**})\overset{\Psi}\to\to E^*\wh\ot_p F
    \overset{j}\to\to N_p(E, F),
$$
where $ j$ is the canonical quotient map. Its kernel
$ \operatorname{ Ker}\,\Psi_0$ consists of all the operator
$ U\in N_p(E, Z^{**})$ which are turned into the identical zero
after the factorization $ \psi: Z^{**}\to F,$ that is of those $ U,$
for which $ U(E)\subset \operatorname{ Ker}\,\psi=\pi_Z(Z).$
This means that
$ \operatorname{ Ker}\, \Psi_0= \NR[p](E, \pi_Z(Z)).$
$\quad\blacksquare$\enddemo

Let us note a curious ``auxiliary" effect of our arguments (this is,
certainly, something like ``to break a butterfly on the wheel"):

\proclaim {\bf Corollary $1'$ {\rm (A. Grothendieck [3])}}\it
If $ p\in (0,2/3]$ then $ G^*\wh\ot_p Y = N_p(G,Y)$ for every
Banach spaces $ G$ and $ Y$ {\rm(}that is there is no factorization under
the canonical quotient map $ G^*\wh\ot_p Y \to N_p(G,Y)${\rm)}.
\endproclaim\rm
\demo{\it Proof}
Clearly, we can consider only separable Banach spaces
\footnote{For each separable $H\subset G^*$ there exists
   a separable subspace $ G_0\subset G$ such that for the identity
   imbedding $j: G_0\to G$ the operator $j^*\big|_H\to G_0^*$
   is an isometric imbedding.}.
If $ p\le 2/3$ then every operator $ U\in \NR[p](G, \pi_Z(Z)),$
where $ Z$ is from lemma 2 (or, that is the same, from corollary 1),
factors through a nuclear diagonal operator (see below remark 2)
and, therefore, is nuclear itself. Thus
$ \NR[p]\(G,\pi_Z(Z)\)\cap \( N_1(G, \pi_Z(Z))\)= \NR[p]\(G,\pi_Z(Z)\)$
and the corollary 1 is applied.
$\quad\blacksquare$\enddemo

Recall that a Banach space $ X$ has the approximation property
of the order $p, \, p\in (0,+\infty]\,$ (briefly,
$X\in \operatornamewithlimits{AP_p}$), if for every Banach space
$ Z$ the canonical mapping
$Z^*\widehat \otimes _p \,X \to \operatorname{L} \left(Z,X \right)$
is one-to-one (see, e.g., [8], [10]). It is convenient to introduce
some more additional notation and definition connected with the
approximation properties.
Namely, for a pair of the spaces $ X,Y$  and a number $ s\in(0,1)$
we denote by $Y^*\widehat \otimes _{s,\infty} \,X $
the linear space consisting of tensors
$ z\in Y^*\widehat \otimes _1 \,X$ (the projective tensor product),
admitted representations of the kind
$$ z=\sum_{k=1}^\infty \la_k \,y'_k\otimes x_k  \
\text{ where }\ \la_k\searrow,\,\la_k^s=o(1/k),\, \|y'_k\|\,\|x_k\|=1.
$$

\definition {Definition {\rm 1}}
The Banach space $X$ has the property
$\operatornamewithlimits{AP_{s,\infty}}\,$, where $s\in (0,1),\,$
if for each Banach space $ Z$ the canonical mapping
$Z^*\widehat \otimes _{s,\infty} \,X  \to
\operatorname{L} \left(Z,X \right)$
is one-to-one (in other words, $Z^*\widehat\otimes _{s,\infty} \,X =
\operatornamewithlimits{N_{s,\infty}}\left(Z,X \right);$  it must be
clear what the operator space on the right is).
\enddefinition

Note that a space $ X$ has the property
$ \operatorname{ AP}_{s,\infty}$ iff
the canonical mapping
$X^*\widehat \otimes _{s,\infty} \,X  \to
\operatorname{L} \left(X,X \right)$ is one-to-one;
in this case
$X^*\widehat \otimes _{s,\infty} \,X =
\operatornamewithlimits{N_{s,\infty}}\left(X,X \right).$

\proclaim {Theorem {\rm 1}}
\roster
\runinitem "{1)}"
Let $\alpha\in (0,1/2],\, C>0,\,$ and $X$ be a Banach space.
If every finite-dimensional subspace $E$ of the space $X$ is
$C\,( \operatorname{ dim}E^\alpha)$-complemented in $X$ then
$X\in \operatornamewithlimits{AP_{s,\infty}},$
where $1/s=1+\alpha$. In particular {\rm (take $\,\alpha =1/2$),}
every Banach space has the property
$\operatornamewithlimits{AP_{2/3,\infty}},$ and if $X$ is
a subspace or a factor space of some space
$\operatornamewithlimits{L^p}(\mu),$ where $p\in (1,+\infty)$, $p\ne 2$,
then $X\in\operatornamewithlimits{AP_{s,\infty}}$ \  $(1/s=1+|1/p-1/2|).$
\endroster
\roster
\runinitem "{\phantom{vv}2)}"
For each $s\in [2/3,1)$ there exists a separable reflexive space
$Y$ such that $Y\in \operatornamewithlimits{AP_{s,\infty}}$, but
$Y\notin \operatornamewithlimits{AP_{r}}$ for any $r\in (s,1]$.
\endroster
\endproclaim

The {\it Proof} of the part 2) of the theorem is completely analogous to
the proof of the corresponding fact from
[10] (namely, the assertions
1) and 2) of the theorem 5.4 in [10]); the proof of this theorem 5.4 is
literally carried over to our case, but instead of the lemma 5.1 from [10]
one has to use such a strengthening of the part 1 of our theorem 1:

\proclaim {Proposition {\rm 1}}
Let $\alpha\in (0,1/2],\, C>0 ,\,$ and $X$ be a Banach space.
Suppose that for every finite-dimensional subspace $E$ of the space $X$
one can find a finite-dimensional subspace
$F,\, E\subset F\subset X,\, $ such that
$F$ is $\, C\left(\operatorname{dim}E \right)^\alpha$-complemented
in $X$. Then $X\in\operatornamewithlimits{AP_{s,\infty}}$, where
$1/s=1+\a$.
\endproclaim
Again, it is enough to refer to the paper [10].
Our Proposition 1 is proved implicitly in [10], Lemma 5.1:
if we look carefully at the proof of the lemma in [10], then
it is not hard to see that the main thing used in that proof is
the belonging of the tensor $ z$ to the projective tensor product
(in this paper all the tensors have this property)
and the estimation $ \lambda_j^s=o(j^{-1})$ for the coefficients of
$ z$ from the proof ot the mentioned Lemma 5.1.
So, we will not repeat all those arguments from [10] which were used
in the proof of this lemma, and refer the reader to [10].

\proclaim {Theorem {\rm 2}}
\roster
\runinitem "{1)}"
\,Let $\, s\in \(0,1\).
$ If $\, X^*\in \,\ap[s\infty]\ $ or $\, Y^{***}\in \,\ap[s\infty],\ $
then $\nr[s\infty]\left(X,Y\right)\subset \n1\(X,Y\).$
In particular, for every Banach spaces $X\,$ and $Y\,$ we have:
$\nr[{2/3}\infty]\left(X,Y\right)\subset \n1\left(X,Y\right).$
\endroster
\roster
\runinitem "{\phantom{vv}2)}"
There exists a Banach space $Z\,$ with the following properties:\newline
      \phantom{vvvv} {\rm(}i{\rm)}
$Z^{**}\,$ has a basis \ $($therefore, $Z^{**}\in \A)$;\newline
      \phantom{vvvv}  {\rm(}ii{\rm)}
for each $\ s\in \(2/3,1\]$ \quad
$\NR[s] \left(Z^{**},Z\right)\not\subset \n1\(Z^{**},Z\).$
\newline
Thus, the approximation conditions imposed on $X\,$ and $Y\,$
in $1)$ are essential.
\endroster
\endproclaim

\demo{Proof}
   1) Suppose there exists such an operator $ T\in L(X,Y)$ that
$ T\notin N_1(X,Y)$ but $ \pi_Y\,T\in N_{s,\infty}(X,Y^{**}).$
Since either $ X^*$ or $ Y^{**}$  has the property $ \ap[s\infty],$
we have:
$$ N_{s,\infty}(X,Y^{**})=X^*\widehat\otimes_{s,\infty} Y^{**}.$$
Hence the operator $ \pi_Y\,T$ can be identified with a tensor element
$ t\in X^*\widehat\otimes_{s,\infty} Y^{**}\subset
X^*\widehat\otimes_1 Y^{**};$ in addition, by the choice of $ T,$ \
$  t\notin X^*\widehat\otimes_1 Y$ \ (as usually,
$  X^*\widehat\otimes_1 Y$ is considered as a subspace of the space
$  X^*\widehat\otimes_1 Y^{**}$). Therefore, there exists such an operator
$ U\in L(Y^{**},X^{**})=\( X^*\widehat\otimes_1 Y^{**}\)^*$ that
$ \tr U\circ t=\tr \(t^*\circ \( U^*|_{X^*}\) \)=1$ and
$ \tr U\circ \pi_Y\circ z=0$ for each $ z\in X^*\widehat\otimes_1 Y.$
It follows now that $ U\pi_Y=0$ and $ \pi_Y^*\,U^*|_{X^*}=0.$
Indeed, if $ x'\in X^*$ and $ y\in Y$ then
$$ <U\pi_Y\,y,x'> = <y, \pi_Y^*\,U^*|_{X^*}x'>
       = \tr \,U\circ (x'\otimes \pi_Y(y))=0.
$$
Evidently the tensor element $ U\circ t$ induces the operator
$ U\pi_Y T,$ which is identically equal to zero.

If $ X^*\in\ap[s\infty]$ then
 $ X^*\widehat\otimes_{s,\infty} Y^{**}= N_{s,\infty}(X,Y^{**}),$
hence this tensor element is zero what is contradicted to the equality
$ \tr\, U\circ t=1.$

Let now $ Y^{***}\in \ap[s\infty].$ In this case the operator
$$ V:= \( U^*|_{X^*}\)\circ T^*\circ \pi_Y^*: \
       Y^{***}\to Y^* \to X^*\to Y^{***}
$$
uniquely defines a tensor element
$ t_0\in Y^{****}\widehat{\otimes}_{s,\infty} Y^{***}.$
Take any representation $ t=\sum x'_n\otimes y''_n$ for $ t$
as an element of the space $ X^*\widehat{\otimes}_{s,\infty} Y^{**}.$
Denoting for the brevity the operator $ U^*|_{X^*}$ by $ U_*,$
we get:
$$\multline
   Vy'''=U_*\, \( T^*\pi_Y^*\,y'''\) =
    U_*\, \( (T^*\pi_Y^*\pi_{Y^*})\,\pi_Y^*\, y'''\) =
    U_*\, \( (\pi_Y T)^*\,\pi_{Y^*})\,\pi_Y^*\, y'''\) = \\ =
    U_*\, \( (\sum y''_n\otimes x'_n)\,\pi_{Y^*}\,\pi_Y^*\, y'''\)
    =U_*\, \( \sum <y''_n, \pi_Y^*\, y'''> \,x'_n \) = \\ =
      \sum <\pi_Y^{**}y''_n,  y'''> \,U_* x'_n.
  \endmultline
$$

So, the operator $ V$ (or the element $ t_0$) has in the space
$ Y^{****}\widehat\otimes_{s,\infty} Y^{***}$ the representation
$$ V= \sum \pi_Y^{**}(y''_n)\otimes U_* (x'_n).
$$
Therefore
$$  \tr t_0=\tr V= \sum <\pi_Y^{**}(y''_n), U_* (x'_n)> =
        \sum <y''_n, \pi_Y^*\,U_* x'_n> =  \sum 0=0.
$$
On the other hand,
$$  Vy'''= U_* \( \pi_Y T\)^* y'''= U_*\circ t^* (y''')=
     U_*\, \( \sum <y''_n, y'''> \, x'_n\)=
      \sum <y''_n, y'''> \, U_* x'_n,
$$
whence $ V=\sum y''_n\otimes U_*(x'_n).$   Therefore
$$ \tr t_0=\tr V= \sum <y''_n, U_* x'_n> = \sum <Uy''_n, x'_n>
= \tr U\circ t=1.$$
The obtained contradiction completes the proof of the first part of 1).
The case $ s=3/2$ is reduced to the just proved assertion by using
the first part of Theorem 1.

  2) It is well known that there exist a separable reflexive Banach space
$ X$ and a nonzero tensor element $ z\in X^*\wh\ot_1 X$ such that
$ z\in X^*\wh\ot_s X$ for all $ s>2/3,$  but the assotiated operator
$ \wt z$ is identically zero (see, e.g., [1], [7]).
Thus for each $ s\in (0,1]$ \ $ X^*\wh\ot_s X\neq N_s(X,X).$
By Corollary 1 there are a separable $ Z$  and a homomorphism
$ \ffi: Z^{**}\to X$ such that $ Z^{**}$ has a basis,
$ \operatorname{ Ker}\,\ffi=\pi_Z(Z),\,$ $ \ffi^*(X^*)$ is complemented
in $ Z^{***}$ and $ \NR[s](X,Z)\not\subset N_1(X,Z)$ for all $ s>2/3.$
Since $ \ffi^*(X^*)$ is complemented in $ Z^{***},$ then
$ \NR[s](Z^{**}, Z) \not\subset Z^{***}\wh\ot_1 Z = N_1(Z^{**}, Z)$
for each $ s>2/3.$
 $\Q$ \enddemo

\remark {\bf Remark {\rm 2}}
It is rather easy to show that $\ \NR[2/3]\subset \n1. $
Indeed, if an operator $ T$ acts from a space $ E$
to a space $ Z$ and is $ 2/3$-nuclear as an operator from $ E$
into $ Z^{**}$ then, by the very definition of the $ 2/3$-nuclear operator,
``splitting" the coefficients of its tensor--series expansion into pairs
of appropriately chosen factors, we can expand the operator
$ T$ itself into a product as follows:
$$
 \gather
T:\  E   @>\ \, A\ \,>>  c_0  @>\ \Delta_1\ >>   l^1  @>\ \,j\ \,>>
                l^2  @>\ \Delta_2\ >>  l^1  @>\ B\ >>  Z^{**} \\
     \endgather
$$
where $ \Delta_j$ are diagonal operators (thus,
$ \Delta_1$ is nuclear) and other operators are continuous.
Taking an orthogonal projection in $ l^2$ onto the closure of the
image of the operator $ j\,\Delta_1\,A,$ and considering
$ B$ only on this subspace, we get the same operator
$ T:E\to Z;$  this simple trick shows that
$ T$ is nuclear from $ E$ into $ Z.$
At this moment I have no idea
how to obtain the inclusion $\,\nr[{2/3}{\infty}]\subset \n1\,$
directly by an analogous way \,(maybe, one have to get a factorization
through $\,l_2\,$ not for whole operator but for its finite dimensional
parts by using the fact that every $\,n$-dimensional subspace of a
Banach space is $\,n^{1/2}$-complemented).
\endremark

\proclaim {Corollary {\rm 2}}
Let $X$ be such a reflexive Banach space that
each its $n$-dimensional subspace $\,Cn^\alpha$-complemented in $X$
\,{\rm (}here $C>0, \alpha \in (0,1/2) \,$ are some constances,
$\,n=1,2,\dotso).$
Then for every Banach space $Y$
 $$\nr[s\infty]\left(X,Y\right)\subset \n1\(X,Y\).$$
where $\,1/s=1+\alpha$.
\par In particular, if $\,X\subset L^p(\mu), 1<p<+\infty, p\ne 2, $
\, É $\,1/s=1+|1/2-1/p|$, \, then
 $$\nr[s\infty]\left(X,Y\right)\subset \n1\(X,Y\) \qquad
 \text{for every }Y.$$
\endproclaim

  {\it To prove}\, the corollary it is enough to apply Theorem 1,1) and
Theorem 2,1)  and to use the fact that any $n$-dimensional subspace of
  $\ L^p(\mu) $ \ $Cn^{|1/2-1/p|}\,$-complemented in  $\, L^p(\mu) $.

  \remark {\bf Remark {\rm 3}}
 In Theorem 2,1) and Corollary 2 we obtained even the following inclusion:
  $\nr[s\infty]\left(X,Y\right)\subset \,X^*\widehat \otimes _1 \,Y.$
One has to understand this in the following manner: every operator
$\,T\in \nr[s\infty]\,\(X,Y\)$ is uniquely generated by a tensor element
    from  the space
  $X^*\widehat \otimes _{s,\infty} \,Y^{**}\subset
  X^*\widehat \otimes _1 \,Y^{**}$,
which (the tensor element), in turn, belongs to the space
$X^*\widehat \otimes _1 \,Y$,
a subspace of $X^*\widehat \otimes _1 \,Y^{**}$; so if
$$\alignat 2
&\sum x'_n \otimes y^{\prime\prime}_n &\quad
  &\text{is a representation of
               $T$ in the space $X^*\widehat \otimes _{s,\infty}Y^{**}$},\\
&\sum \tilde x'_n \otimes \tilde y^{\prime\prime}_n &\quad
  &\text{-- \ in the space
               $X^*\widehat \otimes _1 \,Y^{**}$ \ and}  \\
&\sum \Tilde{\Tilde x}'_n \otimes \Tilde{\Tilde y}_n &\quad
  &\text{-- \ in the space
               $X^*\widehat \otimes _1 \,Y$},
\endalignat $$
then for each $\,U\in \operatorname{L} \(Y^{**},X^{**}\)$
$$\tr U\circ T = \sum \langle x'_n, Uy''_n\rangle =
 \sum \langle\tilde x'_n ,U\tilde y''_n \rangle =
 \sum  \langle\Tilde{\Tilde x}'_n,U\Tilde{\Tilde y}_n\rangle $$
(let us note that the trace is well defined since in Theorem 2,1) and
Corollary 2 the natural map
$\nr[s\infty]\,\(X,Y^{**}\) \to \operatorname{L} \(X,Y^{**}\)$
is one-to-one).
\endremark
By the same method as in the proof of Theorem 2,1) we get
\proclaim {Theorem {\rm 2$'$}}
Let $s\in \(0,1\].$ If $X^*\in \AP s,$ or $Y^{***}\,\in\AP s$  then
$\NR[s] \(X,Y\) \subset \n1\(X,Y\).$
\endproclaim
Looking at the proof of the part 2) of Theorem 2 and of Corollary 2 it is
not hard to see that we have the following stronger result:
\proclaim {Proposition {\rm 2}}
There exist a separable Banach space $Z$ and an operator
$A\in \operatorname{L} \(Z^{**},Z\)$ so that $Z^{**}$ has a basis,
$\pi\, A\in \n s \(Z^{**},Z^{**}\)$ for all $s>2/3$  but
$A\notin \n1 \(Z^{**},Z\).$
\endproclaim

\proclaim {Theorem {\rm 3}}
For each $s\in \[2/3,1\) $ there exists a separable Banach space $Z$
so that
\roster
\item "{(i)}"
$Z^{**}$ has a basis;
\item "{(ii)}"
all the spaces $Z, Z^*, Z^{**},\dotso$ have the property $\AP s$;
\item "{(iii)}"
$Z^{***}$ does not have the property $\AP r$ if $r\in \(s,1\]$;
\item "{(iv)}"
for every Banach space $X$
 $$\nr[s\infty] \(X,Z\) \subset \n1\(X,Z\) \quad \text{and}\quad
     \nr[s\infty] \(Z^{**},X\) \subset \n1\(Z^{**},X\);$$
\item "{(v)}"
 for each $\ r\in \(s,1\]$ \quad
 $\NR[r] \(Z^{**},Z\) \not\subset \n1\(Z^{**},Z\).$
Moreover, there exists an operator
              $U: Z^{**} \to Z $
such that $U\in \NR[r]\(Z^{**},Z \)$  for all $\ r\in \(s,1\],$
but $U\notin \n1\(Z^{**},Z\).$
\endroster
\endproclaim

The method of the proof of Theorem 3 is essentially the same as one
for the part 2) of Theorem 2 by using Theorem 1,2), Corollary 1 and
Theorem 2,1). We omit the proof leaving it to the reader as an not
difficult exercise.
\smallpagebreak

Let us look now at the situation when there are no approximation
restrictions on Banach spaces under consideration  ( as, for instance,
in the case $s=2/3$ in Theorem 2,1) and pose the question of what assertion
can be obtained instead of one in the part 1) of Theorem 2.

The following fact we shall use is proved essentially in [10]
(Theorem 2,1,A).
It is convenient to denote by $\varphi ^{\infty}_{sp}$ the canonical map
$X^*\widehat \otimes _{s,\infty} \,Y  \to  X^*\widehat \otimes _p \,Y$
and by $j_p$ the canonical map $X^*\widehat \otimes _p \,Y \to
\operatorname{L} \(X,Y\) $.

\proclaim {Proposition {\rm 3}}
If $s\le 1$ and $p\ge 2s/(2-s)$ then the mapping
$$i_p:\, \varphi ^{\infty}_{sp}\(X^*\widehat \otimes _{s,\infty} \,Y\) \to
\operatorname{L} \(X,Y\) $$
is one-to-one.
\endproclaim

One can find a proof of this proposition in [10] (see there the proof of the
part $ A$ of Theorem 2,1).
\smallpagebreak

The following fact is a stronger version of the theorem 3.1),A) from [10]
where it was shown that if $s$ and $p$ are as in Proposition 3 then
$\NR[s] \(X,Y\) \subset \n p\(X,Y\) $. However, the proof will be not so
elementary as in [10].

\proclaim {Theorem {\rm 4}}
If $2/3\le s<1$ and $p=2s/(2-s)$ then
$$\nr[s\infty] \(X,Y\)  \subset \n p\(X,Y\) $$
for any Banach spaces $X$  and $Y$.
\endproclaim

\demo{Proof} Let  $ T\in \nr[s\infty]$   and $ T\neq 0.$
Fix a tensor element
$ t=\sum \bar x'_n\ot \bar y''_n\in X^*\wh\ot_{s,\infty} Y^{**},$ which
generates the operator $ T$ (i.e., $\pi_Y T=\wt t$).
It is enough to prove that $ t\in X^*\wh\ot_p Y$ (i.e, that the image
$ \ffi^\infty_{sp}(t)\in X^*\wh\ot_p Y^{**}$ lies in fact
in the subspace $ X^*\wh\ot_p Y$).

Suppose that $ t\notin  X^*\wh\ot_p Y.$  In this case one can find
an operator
$ U\in \Pi_{p'} (Y^{**}, X^{**}),$ for which $ \tr U\circ t=0$ and,
on the other hand, $ \tr U\circ z=1$ for all $ z\in X^*\wt\ot_p Y.$
Further on it will be convenient to consider several diagrams (with some
comments):
$$
 \gather
  X   @>\ \, T\ \,>>   Y   @>\   \pi_Y  >>  Y^{**} @>\ \,B\ \,>>
                Z_{p'}  @>\ S\ >>  X^{**}  @>\ T^{**}\ >>  \pi_Y(Y) \\
     \endgather
$$
($U=SB;\ \wt t=\pi_Y T; \ B\pi_Y=jB_0$\, and
$A:=\pi^{-1}T^{**}S $ (these operators will appear later);
\  $ U\pi_Y=0 $ and $Sj=0);$
$$
 \gather                   Y  @>\ B_0\ >>
  Y_{p'}  @>\ \, j\ \,>>  Z_{p'}  @>\ S\ >>   X^{**}  @>\ \,T^{**}\ \,>>
                \pi_Y(Y)  @>\ \pi^{-1}_Y\ >>  Y  @>\ B_0\ >>  Y_{p'}. \\
     \endgather
$$
Here $ Z_{p'}\subset L_{p'}(\mu)$ (for some finite measure $ \mu$);
$ B\in\Pi_{p'}(Y^{**}, Z_{p'})$ and $S\in L(Z_{p'}, X^{**})\,$ are
such that $ U=SB;$
$Y_{p'}=B\pi_Y(Y);$\ $ B_0$ is the operator induced by $ B;$
$ j$ is the identity imbedding.
Note that $T^{**}(X^{**})\subset \pi_Y(Y)$
(by the compactness of $ T).$
Besides, $ U\pi=0$ by the assumption; hence, $ Sj=0.$
Denote the operator $ \pi^{-1}_Y T^{**}S$ by $ A.$

Since $ Z_{p'}\subset L_{p'}(\mu)$  and
$ \pi_Y A\in N_{s,\infty},$ then, by Theorem 1,1), the operator $ \pi_Y A$
is uniquely defined by a tensor element from
$ Z^*_{p'}\wh\ot_{s,\infty}Y^{**},$
which will be denoted again by $ \pi_Y A.$ It is clear that
$$ \pi_Y A= T^{**} S = \sum S^* \( \pi_{X^*}(\bar x'_n)\) \ot \bar y''_n.
$$
By Corollary 2 and Remark 3, $ A\in Z^*_{p'}\wh\ot_1 Y,$
and also
$$ \multline
 \tr B(\pi_Y A)=\tr (B\pi_Y)\circ A =
\sum <S^*\pi_{X^*}\bar x'_n,\, B\bar y''_n> =\\
=\sum <\bar x'_n,\, U\bar y''_n> = \tr U\circ t=1.
 \endmultline
$$
In particular, it follows from here that $ \pi_Y A\neq 0$\, (and $ A\neq 0$).

Take a representation of $ A$ of the kind $ \sum z'_n\ot y_n$ as
an element of the space
$ Z^*_{p'}\ot_1 Y;$ this is also a representation of the tensor element
$ \pi_Y A\in Z^*_{p'}\wh\ot_{s,\infty} Y^{**}$\ (but in the space
$ \pi_Y A\in Z^*_{p'}\wh\ot_1 Y^{**} $). We have:
$$ \multline
   1= \tr (B\pi_Y)\circ A=\tr (jB_0)\circ A= \tr (jB_0)\circ
\( \sum z'_n\ot y_n\)=\\ = \sum <z'_n,\, jB_0 y_n > =
\sum <z'_n,\, jB_0^{**}\pi_y y_n> = \tr \( jB_0^{**}\)\circ
\( \sum z'_n\ot \pi_Y y_n\) =\\ =
         \tr \( jB_0^{**}\)\circ \( \pi_Y A\)
= \tr \( jB_0^{**}\)\circ
\( \sum S^* \( \pi_{X^*}(\bar x'_n)\) \ot \bar y''_n\)= \\=
\sum <S^*\pi_{X^*}\bar x'_n,\, jB_0^{**}\bar y''_n> =
\tr B_0^{**}\circ \( \sum \((Sj)^*\, \pi_{X^*}\bar x'_n\)\ot \bar y''_n\).
  \endmultline
$$
It follows from here that the tensor element
$ \a:= \sum \((Sj)^*\, \pi_{X^*}\bar x'_n\)\ot \bar y''_n,$
belonging to the space $ Y^*_{p'}\wh\ot_{s,\infty}  Y^{**},$
is not equal to zero.
By Theorem 1,1), the associated operator $ \wt \a$ is not zero too.
It is remain to note that $ \wt\a=T^{**}Sj=0.$
$\Q$ \enddemo

\remark {\bf Remark {\rm 4}}
The assertion $\text{A}'$ of Theorem 3.1 from [10] shows that Theorem 4
is sharp. An analogous remark can be made about Proposition 3. However
the following result is much stronger than the part $\text{A}'$
of the theorem 3.1 in [10] where it is obtained that there are
two Banach spaces $X$ and $Y$ such that, in particular,
$\NR[s] \(X,Y\)\not\subset \n p \(X,Y\)$ and $Y\in \A$.
Having in mind the result of A. Grothendieck (which we said about
before the definition 1 and which is now clear to be partially wrong --
see [6]) I wrote in [10] that the space $Y$  "certainly does not have the
$\operatorname{BAP}$". Now we shall see that this last phrase
was not reasonable.
    \endremark

\proclaim {Theorem {\rm 5}}
Let $r\in [2/3,1)$ and $1/r=1/2+1/p.$
There exist two separable Banach spaces
$Z$ and $W$ such that
\roster
\item "{(i)}"
$W$ and $Z^{**}$ have bases;
\item "{(ii)}"
all the duals of the space $W$ have the property $\ap[r\infty]$;
\item "{(iii)}"
$W^*$ does not have $\AP s$ for any $s \in (r,1]$;
\item "{(iv)}"
for each Banach space $\ E$ \quad $\nr[r\infty]\(W,E\)\subset \n1\(W,E\)$;
\item "{(v)}"
$\NR[s]\(W,Z\)\not\subset \n p \(W,Z\)$ for any $s\in (r,1]$.
\endroster
\endproclaim

\demo{Proof}
It is enough to prove that for every $ s, r<s<1,$ there are Banach spaces
$ Z=Z_s$  and $ W=W_s$ with the properties mentioned in the
statement of the theorem; moreover, if
$ s$ is fixed, it is enough to get in (iii) the weaker condition
$ W^*\notin \ap[s\infty]$ (since $ s$ is running over the interval
$ (r,1)$, and the property $ \AP 1$ is the strongest one among all of the
approximation properties). If we do this then for $ s_n\searrow r$
the spaces $ \( \sum Z_{s_n}\)_{l_2}$ and  $ \( \sum W_{s_n}\)_{l_2}$
will be completely satisfied all the conditions of the theorem.

Thus, let $ s\in (r,1),$ so $ p<2s/(2-s).$
By Theorem 2.1,A'),1) in [10], there exist separable reflexive
spaces $ X$ and $ Y,$  a tensor element $ z\in X^*\wh\ot_s Y$
and an operator $ U\in QN_{p'}(Y,X)$ (quasi-$p'$-nuclear) such that
 $ \tr U\circ z=1$ and the associated operator $ \wt z=0;$ moreover,
the space $ X$ can be taken as a subspace of a space
$ L_{p'}$ (the last fact is contained explicitly in the proof
of the mentioned theorem from [10]).

Let $ Z$ be such a separable space that $ Z^{**}$  has a basis and
there is a homomorphism $ \ffi$ from $ Z^{**}$ onto $ Y$ with the kernel
$ \pi_Z(Z)$ (see Lemma 2). Lift up the (nuclear) tensor element
$ z,$ lying in $ X^*\wh\ot_s Y,$
to get an element $ \a\in X^*\wh\ot_s Z^{**},$
so that $ \ffi\circ \a=z,$ and set $ V:= U\circ \ffi.$
Since $ \tr V\circ\a=\tr U\circ z=1$ and $ Z^{**}$ possesses the
approximation property (the property $ \AP 1$), then $ \wt \a=\a\neq 0.$
Besides, the operator $ \wt{\ffi\circ\a},$ associated with the tensor
$ \ffi\circ\a,$ is equal to zero. Therefore
$ \a(X)\subset \operatorname{ Ker}\ffi= \pi_Z(Z),$ that is the operator
$ \a$ acts from $ X$ to $ Z.$

If $ \a\in X^*\wh\ot_p Z,$ then for any its $ p$-nuclear representation
$ \a=\sum x'_n\ot z_n$ we have:
$ \tr V\circ\a=\sum <x'_n,\, Vz_n>= \sum 0=0;$ hence
$ \a\notin X^*\wh\ot_p Z$ (we used the fact that
$ V\in\Pi_{p'}(Z,X)$). Thus
$\NR[s]\(X,Z\)\not\subset \n p \(X,Z\).$

Let us note that $ X$  (and $ X^*$) has the property $ \ap[r\infty]$
(Theorem 1,1) --- remember where the space $ X$ lies),
but does not have the property $ \ap[s\infty]$ (Theorem 2,1)).

Now, consider such a separable Banach space $ E$ that
$ E^*$ has a basis and there is a homomorphism $ \psi$ from $E^*$
onto $ X$ with the property that $ E^{**}=\psi^*(X^*)\oplus \pi_E(E)$
(Lemma 1). Set $ W=E^*.$ Since $ E\in\AP 1$ and
$ X^*\in \ap[r\infty],$ all the conjugate spaces
$W^*, W^{**},\dots$ possess the property $ \ap[r\infty];$
by Theorem 2,1), the assertion (iv) is fulfilled, and by the construction
the same is true for (v). Besides, $ W^*\notin \ap[s\infty].$
$\Q$ \enddemo

It follows from (i) and (v) of the previous theorem

\proclaim {\bf Corollary 3}\it
There exists a pair of separable Banach spaces $ Z, W$
with the following properties.
The spaces $ Z^{**}$ and $ W$ have bases {\rm (}therefore,
have the approximation
property{\rm )}, and for every $ p, 1\le p<2,\,$
one can find non-$p$-nuclear operator from $ W$ to $ Z$
with $ p$-nuclear second adjoint.
\endproclaim\rm

\remark {\bf Remark 5}
Such a phenomenon is impossible if either
$ W^*$  or
$ Z^{***}$ has the approximation property. This was mentioned above
in passing
(a particular case was formulated in Theorem $2'$).
As to the case $ p>2,$ it is the theme of the forthcoming paper of the
author.
\endremark\medpagebreak
    \vskip0.2in

\centerline{REFERENCES}
\bigpagebreak

\ref \no 1\by A. M. Davie \pages  261-266
\paper  The approximation problem for Banach spaces
\yr 1973\vol 5
\jour Bull. London Math. Soc.
\endref

\ref \no 2\by T. Figiel, W. B. Johnson \pages 197-200
\paper  The approximation property does not imply  the bounded
   approximation property
\yr 1973\vol 41
\jour   Proc. Amer. Math. Soc.
\endref

\ref \no3\by A. Grothendieck \pages 196-140
\paper  Produits tensoriels topologiques et espaces nucl\'eaires
\yr 1955\vol  16
\jour  Mem. Amer. Math. Soc.
\endref

\ref \no 4\by J. Lindenstrauss\pages  279-284
\paper  On James' paper ``Separable Conjugate Spaces"
\yr 1971\vol 9
\jour Israel J. Math.
\endref

\ref \no 5\by B. M. Makarov, V. G. Samarskii\pages   122-144
\paper  Weak sequential completeness and close properties of some
  spaces of operators
\yr 1983\vol
\jour in book ``Theory of operators and theory of functions".
Leningrad University
\finalinfo (in Russian)
\endref

\ref \no 6  \by E. Oja, O. I. Reinov  \pages  121-122
\paper  Un contre-exemple \`a une affirmation de A. Grothendieck
\yr  1987 \vol  305
\jour   C. R. Acad. Sc. Paris. --- Serie I
\endref

\ref \no 7 \by A. Pietsch\pages
 \paper  Operator ideals
 \yr 1980\vol
 \jour North-Holland
 \endref

\ref \no 8  \by O. I. Reinov\pages   125-134
\paper   Approximation properties of order p and the existence of
  non-p-nuclear operators with p-nuclear second adjoints
\yr 1982 \vol  109
\jour    Math. Nachr.
\endref

\ref \no 9  \by  O. I. Reinov \pages  115-116
\paper  A simple proof of two theorems of A. Grothendieck
\yr 1983 \vol 7
\jour  Vestnik LGU
\finalinfo (in Russian)
\endref

\ref \no 10  \by O. I. Reinov \pages 145-165
\paper   Disappearing tensor elements in the scale of $ p$-nuclear
operators
\yr 1983\vol
\jour in book ``Theory of operators and theory of functions".
Leningrad University
\finalinfo (in Russian)
\endref

\ref \no 11\by O. I. Reinov  \pages 905-907
\paper Sur les operateurs p-nucl\'eaires entre espaces de Banach avec bases
\yr 1993\vol  316
\jour  C. R. Acad. Sc. Paris. ---  Serie I
\endref

\enddocument